\documentclass{amsart}

\newtheorem{theorem}{Theorem}[section]
\newtheorem{lemma}[theorem]{Lemma}

\theoremstyle{definition}

\theoremstyle{remark}

\numberwithin{equation}{section}

\begin{document}

\title{Enumeration and limit laws of series-parallel graphs
}

\author{Manuel Bodirsky}
\address{Institut f\"{u}r Informatik, Humboldt-Universit\"{a}t zu Berlin, Germany}
\email{bodirsky@informatik.hu-berlin.de}

\author{Omer Gim\'{e}nez}
\address{Departament de Tecnologia, Universitat Pompeu Fabra,
 Passeig de Circumval.laci\'{o}~8,
08003 Barcelona, Spain} \email{Omer.Gimenez@upc.edu}

\thanks{Research supported in part by Projects BFM2001-2340, MTM2004-01728
and Beca Fundaci\'{o} Cr\`{e}dit Andorr\`{a}.}

\author{Mihyun Kang}
\address{Institut f\"{u}r Informatik, Humboldt-Universit\"{a}t zu Berlin, Germany}
\email{kang@informatik.hu-berlin.de}

\author{Marc Noy}
\address{Departament de Matem\`{a}tica Aplicada II, Universitat Polit\`{e}cnica de
Catalunya, Jordi Girona 1--3, 08034 Barcelona, Spain}
\email{Marc.Noy@upc.edu}


\date{\today}


\begin{abstract}
We show that the number $g_n$ of labelled series-parallel graphs
on $n$ vertices is asymptotically  $g_n \sim g\cdot n^{-5/2}
\gamma^n n!$, where $\gamma$ and $g$ are explicit computable
constants. We show that the number of edges in random
series-parallel graphs is asymptotically normal with linear mean
and variance, and that the number of edges is sharply concentrated
around its expected value. Similar results are proved for labelled
outerplanar graphs
 and for graphs not containing $K_{2,3}$ as a minor.
\end{abstract}

\maketitle

\section{Introduction}

A graph is series-parallel (SP for short) if it does not contain
the complete graph $K_4$ as a minor; equivalently, if it does not
contain a subdivision of $K_4$. Since both $K_5$ and $K_{3,3}$
contain a subdivision of $K_4$, by Kuratowski's theorem a SP graph
is planar.
Another characterization, justifying the name, is the following. A
connected graph is SP if it can be obtained from a single edge by
means of the the following two operations: subdividing an edge
(series extension); and duplicating an edge (parallel extension).
In addition, a 2-connected graph is SP if it can be obtained from
a double edge by means of series and parallel extensions; in
particular, this implies that a 2-connected SP graph has always a
vertex of degree two. Although SP operations may give rise to
multiple edges, in this paper all graphs considered are simple.
Yet another characterization is that SP graphs are precisely the
graphs with tree-width at most two. Equivalently they are
subgraphs of 2-trees, where a 2-tree is a graph formed by,
starting from a triangle, adding repeatedly a new vertex and
joining it to an existing edge.

An outerplanar graph is a planar graph that can be embedded in the
plane so that all vertices are in the outer face. They are
characterized as those graphs not containing a minor isomorphic to
(or a subdivision of) either $K_4$ or $K_{2,3}$. They constitute
an important subclass of the class of SP graphs.

Series-parallel graphs have been widely studied in graph theory
and computer science. They are simple in structure but yet rich
enough so that several theoretical and computational problems are
still unsolved on SP graphs. In fact, they are often used as a
benchmark for analyzing the complexity of graph theoretical
problems. The same thing can be said, maybe even more, about
outerplanar graphs.

In this paper we study the enumeration of labelled series-parallel
and outerplanar graphs. From now on, unless stated otherwise, all
graphs are labelled. Next we summarize what is known about this
problem. An SP graph on $n$ vertices has at most $2n-3$ edges.
Those having this number of edges are precisely the 2-trees; it is
known~\cite{moon} that the number of labelled 2-trees on $n$
vertices is equal to ${n \choose 2} (2n-3)^{n-4}$.
 %
On the other hand, an outerplanar graph is 2-connected if and only
it has a unique Hamilton cycle. It follows that a 2-connected
outerplanar graph is in fact equivalent to a dissection of a
convex polygon, the boundary of the polygon being the unique
Hamilton cycle. Hence counting 2-connected outerplanar graphs
amounts essentially to counting dissections of a convex polygon, a
classical and well-known problem that we revisit in Section
\ref{sec:outer}. It is also worth mentioning that an outerplanar
\emph{map} (a map is a planar graph together with a particular
embedding in the plane) on $n$ vertices can be encoded with $3n$
bits \cite{gavoille}; hence the number of outerplanar graphs is at
most $2^{3n} = 8^n$.

The main goal of this paper is to give precise asymptotic
estimates for the number of SP and outerplanar graphs. In
Section~\ref{two-SP} we show that the number $b_n$ of 2-connected
SP graphs is asymptotically of the form
 $$
 b_n \sim b \cdot n^{-5/2} R^{-n} n!,
 $$
where $b$ and $R \approx 0.12800$ are computable constants. Then
in Section~\ref{sec:sp} we show that the total number $g_n$ of SP
graphs is given by
$$
 g_n \sim g \cdot n^{-5/2} \rho^{-n} n!,
 $$
where $\rho \approx 0.11021$.

In Section \ref{sec:edges} we analyze the distribution of  the
number of edges in SP graphs. If $X_n$ is the random variable
equal to the number of edges in a random SP graph, we prove that
$X_n$ is asymptotically normal and that the mean $\mu_n$ and
variance $\sigma_n^2$ of $X_n$  satisfy
$$
\mu_n \sim \kappa n, \qquad \sigma_n^2 \sim \lambda n,
$$
where $\kappa \approx 1.61673$ and $\lambda \approx 0.55347$. As a
consequence, the number of edges is sharply concentrated around
its expected value.

In Section \ref{sec:outer}, we study the same problems for
outerplanar graphs. We prove that the number of outerplanar graphs
is
$$
 h_n \sim h \cdot n^{-3/2} \sigma^{-n} n!,
 $$
where $\sigma \approx 0.13659$, and that the number of edges in
random outerplanar graphs is asymptotically normal with mean
variance
$$
\mu_n \sim \zeta n, \qquad \sigma_n^2 \sim \eta n,
$$
where $\zeta \approx 1.56251$ and $\eta \approx 0.22399$.
Previously, it had been shown \cite{gm} that $\zeta \ge 7/5$.

In Section \ref{sec:comp} we study the distribution of the number
of connected components in series-parallel and outerplanar graphs.
In both cases it turns out that the distribution is asymptotically
a shifted Poisson law with parameter equal to $\nu = 0.117614$ for
SP graphs, and $\xi = 0.14840$ for outerplanar graphs. As a
consequence the probability that a random SP graph is connected
tends to $e^{-\nu} =0.889038$, and to $e^{-\xi}=0.86208$ for
outerplanar graphs.

Finally, in Section \ref{sec:k23} we study the family of graphs
not containing $K_{2,3}$ as a minor. This is a subfamily of
series-parallel graphs which turns out to be very close to the
family of outerplanar graphs.

The proofs are based on singularity analysis of generating
functions and perturbation of singularities (see~\cite{FO,FS}),
and on several ideas developed in \cite{bender} and~\cite{GN} for
solving similar problems for the class of planar graphs. For the
techniques and results of singularity analysis used in the sequel
we refer to the forthcoming book~ \emph{Analytic Combinatorics} by
Flajolet and Sedgewick~\cite{FS}.


\section{Counting two-connected series-parallel
graphs}\label{two-SP}

Let $b_{n,q}$ be the number of 2-connected SP graphs with $n$
vertices and $q$ edges, and let
 $$
    B(x,y) = \sum b_{n,q} y^q {x^n \over n!}
 $$
be the corresponding exponential generating function (EGF).

Following \cite{walsh}, we define a \emph{network} as a graph with
two distinguished vertices, called poles, such that adding the
edge between the poles the resulting multigraph  is 2-connected.
If $D(x,y)$ is the EGF for SP networks, where again $x$ marks
vertices and $y$ marks edges then, as shown in \cite{walsh}, we
have
\begin{equation}\label{eq:B-D}
 {\partial B(x,y) \over \partial y} = {x^2 \over 2}
 \left(
 {1+D(x,y)  \over 1+y}
 \right).
\end{equation}

Since a 2-connected SP graph has always a vertex of degree two, it
follows that there are no 3-connected SP graphs; in the
terminology of \cite{walsh}  there are only s-networks and
p-networks and there are no h-networks. Hence equation (12) in
\cite{bender} simplifies to
\begin{equation}\label{eq:D}
    \log\left( {1+D \over 1+y} \right) = {xD^2 \over 1+xD}.
\end{equation}
Our goal is to perform a complete singularity analysis of $B(x,y)$
using equations (\ref{eq:B-D}) and (\ref{eq:D}). To this end we
first determine the singularities of $D(x,y)$.

From now on $y$ is a fixed positive value. Because of (\ref{eq:D}), the
inverse of $D(x,y)$ as a function of $x$ is given by
 $$
    \psi_y(u) = {\log\left({1+u \over 1+y}\right) \over
    u \left(u -  \log\left({1+u \over 1+y}\right)  \right) }.
 $$
We show that the equation $\psi_y'(u)=0$ has a unique positive
root $u=\upsilon(y)$ for every positive $y$. (We  often write
$\psi$ and $\upsilon$ instead of $\psi_y$ and $\upsilon(y)$ for
brevity.) Hence $D(x,y)$, being the inverse of $\psi$, ceases to
be analytic at $x=R(y)=\psi(\upsilon)$.  By Proposition IV.4 in
\cite{FS}, it follows that the dominant singularity of $D(x,y)$
for fixed $y$ is at $R(y)$. The next result gives a procedure for
obtaining $R(y)$ as a function of $y$. But first we need a
technical lemma, analogous to Lemma~2 in \cite{bender}.

\begin{lemma}\label{one-to-one}
The function
 $$
Y(t) = \frac{1}{1-t^2} \exp\left(\frac{-t^2}{1+t}\right)-1
$$
has an analytic inverse in the domain $t\in(0,1)$.
\end{lemma}
\begin{proof}
It is enough to see that $Y(t)$ is an increasing function in $t\in(0,1)$,
and this is so because
$$
Y'(t) =
\frac{(t+3)t^2}{(t^2-1)^2(1+t)}\exp\left(\frac{-t^2}{1+t}\right)
> 0.
$$
\end{proof}

\begin{theorem}\label{th:singB}
For fixed $y>0$, the dominant singularity of $D(x,y)$
 is at $R(y) = q(t)$, where
\begin{equation}\label{sing}
     q(t)= \frac{(1+t)(t-1)^2}{t^3},
\end{equation}
and $t$ is the unique root of $Y(t) = y$, where $Y(t)$ is as in Lemma
\ref{one-to-one}.
\end{theorem}

\begin{proof}
Let
\begin{equation}\label{def-L}
 L = L(u) = \log\left({1+u \over 1+y}\right).
\end{equation}
A routine computation gives
 $$
 \psi'(u) = {(1+u)L^2 - 2u(1+u)L + u^2 \over (1+u)u^2(L-u)^2}.
 $$
The numerator vanishes when the corresponding quadratic equation
on $L$ has a root, necessarily at $u=\upsilon = \upsilon(y)$.
 This gives
\begin{equation}\label{value-L}
 L(\upsilon) = {2\upsilon(1+\upsilon) - \sqrt{4\upsilon^2(1+\upsilon)^2 - 4\upsilon^2(1+\upsilon)} \over 2(1+\upsilon)}
 = \upsilon-\upsilon^{3/2}(1+\upsilon)^{-1/2}.
\end{equation}
In order to simplify (\ref{value-L}),  we set
$t=\sqrt{\upsilon}/\sqrt{1+\upsilon}$ so that $\upsilon=
t^2/(1-t^2)$.
 Equations (\ref{def-L}) and (\ref{value-L}) then become
 $$
 L(t) = \frac{t^2}{1+t},
 $$
 $$
 L(t) = \log{\frac{1}{(1+y)(1-t^2)}}.
 $$
We solve for $y$ and we observe that $t$ is determined by the
equation $Y(t)=y$. Since the limits of $Y(t)$ when $t$ approaches
$0$ and $1$ are, respectively, $0$ and $+\infty$, it follows from
Lemma~ \ref{one-to-one} that every $y>0$ has a corresponding
$t\in(0,1)$; this determines the unique root $\upsilon$ of
$\psi'(u)$.

The dominant singularity $R(y)$ is at
$\psi(\upsilon)=L/(\upsilon(\upsilon-L))$ which, in terms of $t$,
gives (\ref{sing}).
\end{proof}

Since $\psi'$ has a root $\upsilon$, we know that $D(x,y)$, for
fixed $y$, has a singularity of square-root type. In the following
lemma we show that $\psi''(\upsilon) < 0$ for every $y$, hence the
singular expansion of $D(x,y)$ at the singularity $R(y)$ is (see
Proposition VI.1 in \cite{FS})
 $$
    D(x,y) = \upsilon(y) - \sqrt{ -2R(y) \over \psi''(\upsilon)} X +
    \mathcal{O}(X^2),
    $$
where $X = \sqrt{1-x/R(y)}$.

\begin{lemma}\label{th:exp-D}
For fixed $y>0$, the singular expansion of $D(x,y)$ at
 $R(y)$ is
\begin{equation}\label{eq:exp-D}
    D(x,y) = D_0(y) + D_1(y) X + D_2(y) X +    \mathcal{O}(X^3),
\end{equation}
where $X = \sqrt{1-x/R(y)}$ and
 $$
 D_0 = \frac{t^2}{1-t^2}, \quad  D_1 = \frac{\sqrt{2t^3}}{\sqrt{t+3}(t^2-1)},
\quad D_2 = \frac{2t(t^2+3t+3)}{3(1-t)(3+t)^2}
 $$
where $t$ is the unique root of $Y(t)=y$.
\end{lemma}

\begin{proof}
The constant term is $D_0=D(R(y),y) = \upsilon$, and the terms
$D_1$ and $D_2$ are
 $$D_1 = -
\sqrt{-2 R(y)/ \psi''(\upsilon)}, \quad D_2 =
\frac{R\psi'''(\upsilon)}{3\psi''(\upsilon)^2},
 $$
where these expressions can be obtained by inverting the Taylor
series of $\psi(u)$ at the point $u=\upsilon$ where
$\psi'(\upsilon)$ vanishes.

Now we use Maple to compute $\psi''(u)$ and $\psi'''(u)$, which
are rational expressions on $u$ and $L$. Then we set
$u=t^2/(1-t^2)$ and $L=t^2/(1+t)$, and simplifying we obtain
 $$
 \psi''(t) =  -\frac{(t+3)(t-1)^4(1+t)^3}{t^6},
$$
$$
 \psi'''(t) =  -\frac{2(t^2+3t+3)(t-1)^5(1+t)^5}{t^8}.
 $$
Hence $\psi''(t)<0$ for $t\in(0,1)$, and we use this expressions
to obtain $D_1$ and $D_2$ in terms of $t$.
\end{proof}

After completing the analysis of $D(x,y)$ we turn to that of
$B(x,y)$. The first thing is to express $B$ in terms of $D$.

\begin{lemma}\label{le:B}
The following holds, where $D=D(x,y)$:
\begin{equation}\label{mainB}
   B(x,y) = {1 \over 2} \log(1+xD) - {xD(x^2D^2 +xD + 2 -2x) \over
   4(1+xD)}
\end{equation}
\end{lemma}

\begin{proof}
We follow the proof of Lemma 5 in [GN]. From now on $x$ is a
fixed value. From (\ref{eq:B-D}) it follows that
$$
     B(x,y) = {x^2\over 2}\log(1+y) +
    {x^2 \over 2} \int_0^y {D(x,t) \over 1+t} dt.
$$
Integrating by parts we get
 $$
  \int_0^y {D(x,t) \over 1+t} dt = \log(1+y) D - \int_0^y
  \log(1+t) {\partial D \over \partial t} dt.
 $$

 Now we notice that the inverse of $D$ with respect to $y$ is
 $$
 \phi(u) = -1 + (1+u) \exp\left(-  {xu^2 \over 1+xu} \right).
 $$
The last integral, after the change $s = D(x,t)$ becomes
 $$
 \int_0^{D(x,y)} \left( \log(1+s) - {xs^2 \over 1+xs} \right) ds,
 $$
which can be integrated in elementary terms. The rest of the
computation is routine and the claim follows.
\end{proof}

In view of the expression in Lemma \ref{le:B}, it is clear that,
for fixed $y$, the dominant singularity of $B(x,y)$ is the same as
that of~$D(x,y)$, namely $R(y)$. Using (\ref{mainB}) we can find
the singular expansion of $B(x,y)$.

\begin{lemma}\label{th:exp-B}
For $y$ fixed closed to $1$, the singular expansion of $B(x,y)$ at
its singularity $R(y)$ is
\begin{equation}\label{eq:exp-B}
    B(x,y) = B_0(y) + B_2(y) X^2 +   B_3(y) X^3 +  \mathcal{O}(X^4),
\end{equation}
where $X = \sqrt{1-x/R(y)}$ and $B_0(y),B_2(y),B_3(y)$ are the
following analytic functions of the unique root $t$ of $Y(t)=y$,

$$\begin{array}{rl}
 B_0(t) = & \displaystyle
 \frac{t^3+2\ln(1/t)t^3+2t^2-5t+2}{4t^3} \\[2ex]
 B_2(t) = & \displaystyle
 \frac{(t-1)^3(t+2)}{2t^3} \\[2ex]
 B_3(t) = & \displaystyle
 (1-t)^3\sqrt{\frac{2}{3(t+3)t^3}}
 \end{array}
$$
\end{lemma}

\begin{proof}
It is enough to set $x=R(1-X^2)$ and $D=D_0+D_1X+D_2X^2+D_3X^3$ in
(\ref{mainB}). All the calculations have been performed with
Maple. In particular, we obtain that $B_1$ vanishes identically as
a function of $t$, and that $B_3$ does not depend on the value of
$D_3$.
\end{proof}

\begin{theorem}
The number of 2-connected SP graphs $b_n$ is asymptotically
$$
 b_n \sim b\cdot n^{-5/2}\cdot R^{-n} n!
$$
where $R=R(1)\approx 0.12800$ and $b\approx 0.0010131$.
\end{theorem}
\begin{proof}
By transfer theorems, the asymptotic estimate follows from the
singularity expansion of $B(x,1)$ of Lemma~\ref{th:exp-B}. Solving
$Y(t)=1$ gives $t \approx 0.80703$, and from here we obtain the
values of $R(1)=q(t)$ and of $b=3B_3(t)/(4\sqrt{\pi})$.
\end{proof}


\section{Counting series-parallel graphs}\label{sec:sp}

Recall that  $g_n$, $c_n$ and $b_n$ denote, respectively, the
number of SP graphs, connected SP graphs, and 2-connected SP
graphs on $n$ vertices. Adapting the proof of Lemma 1 in
\cite{GN}, we obtain that  the corresponding exponential
generating functions are related as follows.

\begin{lemma}\label{le:BCG}
The series $G(x)$, $C(x)$ and $B(x)$ satisfy the following
equations:
 $$
    G(x) = \exp(C(x)), \qquad xC'(x) =
    x\exp\left(B'(xC'(x))\right),
 $$
where $C'(x) = {\rm d}C(x)/{\rm d}x$ and $B'(x) = {\rm d}B(x)/{\rm
d}x$.
\end{lemma}

 Let $b_{n,q}$ be the number of 2-connected planar graphs with
$n$ vertices and $q$ edges, and let
 $$
    B(x,y) = \sum b_{n,q} y^q {x^n \over n!}
 $$
be the corresponding bivariate generating function. Notice that
$B(x,1) = B(x)$. The generating functions $C(x,y)$ and $G(x,y)$
are defined analogously. Since the parameter ``number of edges''
is additive under taking connected and 2-connected components, the
previous lemma can be extended as follows.

\begin{lemma}\label{le:BCG-biv}
The series $G(x,y)$, $C(x,y)$ and $B(x,y)$ satisfy the following
equations:
 $$
    G(x,y) = \exp(C(x,y)), \qquad x {\partial \over \partial x } C(x,y)    =
    x\exp\left({\partial \over \partial x} B(x {\partial \over \partial x } C(x,y),y)
    \right).
 $$
\end{lemma}

Let  $F(x,y) = xC'(x,y)$, where the derivative is with respect to
the first variable.  Lemma \ref{le:BCG} implies that
 $$
 F(x,y) = x \exp(B'(F(x,y),y)).
 $$
 It follows that, for fixed $y$, the functional inverse of $F(x,y)$ is
$$
  \Psi_y(u) = u e^{-B'(u,y)}.
$$
The function $ \Psi_y$ should not be confused with $\psi_y$ in
Section~\ref{two-SP}, although it plays a similar role. Our goal
is to prove that for each $y>0$, $\Psi_y'(u)=0$ has a root
$\tau(y)$. As in the previous section we often omit the fact that
$\Psi$ and $\tau$ depend on a fixed~$y$.

\begin{lemma}\label{eq-def-tau}
The equation
\begin{equation}\label{pol-tau}
D^6u^4+D^5u^3+2D^3u^2+4D^2u-2=0,
\end{equation}
where $D=D(u,y)$ and $y$ is a fixed positive value, has a unique
solution $u=\tau(y)$ in $(0,R(y))$.
\end{lemma}
\begin{proof}
Let $T(u,D)=T(u,y)$ be the left hand side of (\ref{eq-def-tau}),
which is an increasing function of $u$ since $D(u,y)$ has
non-negative coefficients. Since $T(0,y)=-2$, it follows that
$\tau$ exists and is unique if and only if $T(R(y),y)>0$. We use
the expressions in terms of $t$ for $R(y)$ and $D(R(y),y)$ given
in Theorem~\ref{th:singB} and Lemma~\ref{th:exp-D} and, after
simplification, we obtain that $T(R(y),y)$ written as a function
of $t$ is
$$
\frac{1-t}{(1+t)^2}.
$$
This is a positive value when $t\in(0,1)$, so the claim follows.
\end{proof}

\begin{theorem}\label{th:singF}
Let $y$ be a fixed positive value. The unique root of $\Psi'(u)=0$
is given by $\tau(y)$ in Lemma~\ref{eq-def-tau}. The dominant
singularity of $F(x,y)$ is at $\rho(y)$, where $\rho$, as a
function of $\tau$, is
\begin{equation}\label{singF}
     \rho(\tau)=\tau\exp\left(\frac{\tau D(\tau D^2-2)}{2(1+\tau D)}\right),
\end{equation}
where $D=D(\tau(y),y)$.

The singular expansion of $F(x,y)$ at its dominant singularity
$\rho(y)$ is
\begin{equation}\label{eq:exp-F}
    F(x,y) = F_0(y) + F_1(y) X + \mathcal{O}(X^2),
\end{equation}
where $X = \sqrt{1-x/\rho(y)}$ and

 $$\begin{array}{rl}
 F_0(\tau) = & \displaystyle \tau \\
 F_1(\tau) = & \displaystyle 2\frac{1-2\tau D^2-\tau^2 D^3}{D}\sqrt{\frac{\tau (1+\tau D)}{S}}\\[2ex]
 S = & -4\tau^5D^7 -5\tau^4D^6 +(6\tau^4-\tau^3)D^5 +5\tau^3D^4
-3\tau^2D^3 +\\
  & 6\tau^2 D^2 +12\tau D +4
 \end{array}
  $$

\end{theorem}

\begin{proof}
We start by differentiating $\Psi(u)$:
$$
\Psi'(u)=\exp(-B'(u,y))(1-uB''(u,y)).
$$
By Lemma~\ref{le:B}, the functions $\Psi(u)$ and $\Psi'(u)$ can be
written in terms of $D=D(u,y)$,
$$
\Psi(u)=u\exp\left(\frac{uD(uD^2-2)}{2(1+uD)}\right)
$$
$$
\Psi'(u)=\frac{u^4D^6+u^3D^5+2u^2D^3+4uD^2-2}{(2u^2D^3+4uD^2-2)(1+uD)}\,
\exp\left(\frac{uD(uD^2-2)}{2(1+uD)}\right  ),
$$
where $D=D(u,y)$. To obtain the previous expressions we have used
the relation
$$
D'={D^2(1+D) \over1-2uD^2-u^2D^3},
$$
which follows directly from (\ref{eq:D}).

Clearly $\Psi'(u)$ vanishes at the roots of the polynomial
$$
  T(u,D)=u^4D^6+u^3D^5+2u^2D^3+4uD^2-2,
$$
hence the root $u=\tau(y)$ of $\Psi'(u)$ is the one given by
Lemma~\ref{eq-def-tau}.

As for the remaining expressions, $\rho(y)$ is $\Psi(\tau)$,
$F_0(y)$ is just $\tau$, and $F_1(y)$ is given by
$-\sqrt{-2\Psi(\tau)/\Psi''(\tau)}$, if we can show that
$\Psi''(\tau)<0$. To obtain $\Psi''(\tau)$ we differentiate
$\Psi'(u)$ with respect to $u$. Note that, when evaluating at
$u=\tau$, the polynomial $T(u,D)$ vanishes, and so

$$ \Psi''(\tau)=\frac{\frac{\partial T}{\partial
u}(\tau,D)+\frac{\partial T}{\partial D}(\tau,D)D'} {(2\tau^2
D^3+4\tau D^2-2)(1+\tau D)}
 \exp\left(\frac{\tau D(\tau
D^2-2)}{2(1+\tau D)}\right).
$$

\smallskip

All factors in this expression are positive but for $2\tau^2
D^3+4\tau D^2-2<T(\tau,D)=0$, hence we have shown that
$\Psi''(\tau)<0$. The expression for $F_1(\tau)$ follows by
straightforward simplification.
\end{proof}

In order to find the singular expansion of $C(x,y)$, we start with
a simple lemma.

\begin{lemma}\label{le:C}
The following holds:
\begin{equation}
   C(x,y) = F(x,y)(1+\log x - \log F(x,y))+B(F(x,y),y)
\end{equation}
\end{lemma}

\begin{proof}
This result is given in the proof of Theorem 1 in [GN]. It
is analogous to that of Lemma~\ref{le:B}, but simpler.
Since $F(x,y)=xC'(x,y)$, we have that
$$
  C(x,y)=\int_0^x \frac{F(s,y)}{s} \, ds=
F(x)\log x-\int_0^x F'(s,y)\log s \, ds.
$$
Now we change variables $t=F(s)$, so that $s=\Psi(t)=t\exp(-B'(t,y))$.
Then the last integral becomes
$$
\int_0^{F(x,y)} \log \Psi(t) \,dt= \int_0^{F(x,y)} (\log t -
B'(t,y)) \,dt.
$$
Hence
$$
C(x,y)=F(x,y)(1+\log x -\log F(x,y))+B(F(x,y),y).
$$
\end{proof}

\begin{theorem}
Let $y$ be a fixed positive value. The dominant singularities
of $C(x,y)$ and $G(x,y)$ are at $\rho(y)$, where $\rho(y)$ is as
in Theorem~\ref{th:singF}. The singular expansions of $C(x,y)$
and $G(x,y)$ at their singularities are
$$\begin{array}{rl}
     C(x,y) = & C_0(y) + C_2(y) X^2 + C_3(y) X^3 + \mathcal{O}(X^4), \\
     \\
     G(x,y) = & G_0(y) + G_2(y) X^2 + G_3(y) X^3 +
     \mathcal{O}(X^4),
\end{array}
$$
where $X = \sqrt{1-x/\rho(y)}$ and
$$
\begin{array}{lll}
 C_0=\tau(\log \rho-\log(\tau)+1)+B(\tau,y), & C_2=-F_0, &
 C_3=-\frac{3}{2}F_1,
\\[1ex]
 G_0=\exp(C_0), & G_2=\exp(C_0)C_2, & G_3=\exp(C_0)C_3.
\end{array}
$$
\end{theorem}

\begin{proof}
It is clear that $G$ and $C$ have the same singularities than $F$.
The singular expansion of $F(x,y)=xC'(x,y)$ can be obtained from
that of $C(x,y)$ by differentiating  and multiplying by
$x=\rho(y)(1-X^2)$, so that
$$
  F(x,y) = (-C_2(y) - \frac{3}{2}C_3(y)X)(1-X^2)+\mathcal{O}(X^2).
$$
By equating coefficients the expressions for $C_2$ and $C_3$
follow. To obtain $C_0$ we evaluate $C(x,y)$ at its singularity
$x=\rho(y)$ using Lemma~\ref{le:C}, and notice that
$F(\rho(y),y)=\tau(y)$.

Finally, the singular expansion of $G(x,y)$ follows from
$G(x,y)=\exp(C(x,y)$, since
$$
\begin{array}{rl}
  G(x,y) & =\exp(C_0)\exp(C_2X^2+C_3X^3)+\mathcal{O}(X^4) \\[1ex]
  & =\exp(C_0)(1+C_2X^2+C_3X^3)+\mathcal{O}(X^4).
\end{array}
$$
\end{proof}

\begin{theorem}
The number of connected SP graphs $c_n$ and all SP graphs $g_n$
are asymptotically
$$
 c_n \sim c\cdot n^{-5/2}\cdot \rho^{-n} n!
$$
$$
 g_n \sim g\cdot n^{-5/2}\cdot \rho^{-n} n!
$$
where $\rho=\rho(1)\approx 0.11021$, $c\approx 0.0067912$ and
$g\approx 0.0076388$.
\end{theorem}

\begin{proof}
The asymptotic estimates follow again from transfer theorems on
the generating functions $C(x,1)$ and $G(x,1)$. As for the
constants, solving Equation~\ref{pol-tau} in
Lemma~\ref{eq-def-tau} for $y=1$ yields $\tau(1)\approx 0.1279695$
and $D(\tau(1),1)\approx 1.84351$, and from here follow the values
of $\rho=\Psi(\tau)$, $c=3C_3/(4\sqrt{\pi})$ and $g=\exp(C_0)\,c$.
\end{proof}

\section{The number of edges in series-parallel graphs}\label{sec:edges}

The main tool in this section is the so-called Quasi-Powers
theorem \cite{FS}, which allows to deduce a normal limit law for a
combinatorial parameter from the bivariate singular expansion of
the corresponding generating function. The proof scheme is as for
Theorem 2 from \cite{GN}. The exact form we need of the
Quasi-Powers theorem is that of Proposition 2 in \cite{GN}.

We work out in some detail the case of 2-connected SP graphs; the
remaining cases follow the same pattern. We know that for fixed
$y$ we have a singular expansion
$$
    B(x,y) = B_0(y) + B_2(y) X^2 +   B_3(y) X^3 +  \mathcal{O}(X^4),
$$
where $X = \sqrt{1-x/R(y)}$ and the $B_i$ are analytic functions.
We deduce that number of edges in 2-connected graphs is normally
distributed and that the expected number of edges is
asymptotically $\alpha n$, where
 $$
\alpha =  -{R'(1) \over R(1)}  \approx 1.71891.
 $$
The derivative $R'(1)$ is computed using the explicit form of
$R(y)$ given in Theorem~\ref{th:singB}; indeed $R'(y) =
q'(t)/Y'(t)$, where $t$ is the unique solution of
 $Y(t) = y$. The relevant values are
$R(1) = 0.12800$ and $R'(1)  = -0.22002$.

The variance is asymptotically
 $\beta n$, where
 $$
 \beta =  -{R''(1) \over R(1)}  -{R'(1) \over R(1)}
    + \left( {R'(1) \over R(1)}\right)^2.
 $$
We compute $R''(1)  = 0.57667$, so that $\beta \approx 0.16846$.
Hence we have proved:

\begin{theorem}\label{th:edges-2con}
Let $X_n$ denote the number of edges in a random 2-connected
series-parallel graph with $n$ vertices. Then $X_n$ is
asymptotically normal and the mean $\mu_n$ and variance
$\sigma_n^2$ satisfy
\begin{equation}
\mu_n \sim \kappa_0 n, \qquad \sigma_n^2 \sim \lambda_0 n,
\end{equation}
where $\kappa_0 \approx 1.71891$ and $\lambda_0 \approx 0.16846$.
\end{theorem}

It is worth recalling that the number of edges in a 2-connected SP
graph is between $n$ and $2n-3$.

For connected SP graphs and arbitrary SP graphs the same result
holds, but in this case the dominant singularity is at $\rho(y)$,
which is given in Theorem~\ref{eq:exp-F}.

\begin{theorem}\label{th:edges}
Let $X_n$ denote the number of edges in a random series-parallel
graph with $n$ vertices. Then $X_n$ is asymptotically normal and
the mean $\mu_n$ and variance $\sigma_n^2$ satisfy
\begin{equation}\label{mean-var}
\mu_n \sim \kappa n, \qquad \sigma_n^2 \sim \lambda n,
\end{equation}
where $\kappa \approx 1.61673$ and $\lambda \approx 0.21125$. The
same is true, with the same constants, for \emph{connected} random
SP graphs.
\end{theorem}

Since $\sigma_n = o(\mu_n)$ it follows that the number of edges in
random SP is concentrated around it expected value, in the sense
that for every $\epsilon >0$ we have
 $$
    \hbox{Prob} \{|X_n - \kappa n| > \epsilon n\} \to 0,
  \qquad \hbox{as $n \to \infty$}.
  $$
This comment also applies to Theorem \ref{th:edges-2con} and to
the Gaussian limit laws presented in the next sections.

\begin{proof}
Since $\rho(y) = \Psi(\tau(y), y)$ it follows that
$$
\rho'(y) = \frac{\partial\Psi}{\partial x}(\tau(y),y)\tau'(y)+
\frac{\partial\Psi}{\partial
y}(\tau(y),y)=\frac{\partial\Psi}{\partial y}(\tau(y),y),
$$
where the first summand vanishes by definition of $\tau(y)$. We
can compute $\partial \Psi/\partial y$ explicitly by
differentiating $\Psi(u,y)$ with respect to $y$, and using that
$$
\frac{\partial D}{\partial y}(x,y)=
-\frac{(1+xD(x,y))^2(1+D(x,y))}{(-1+2xD(x,y)^2+x^2D(x,y)^3)(1+y)}.
$$
We obtain that $\rho'(1)\approx -0.17818$.

To compute $\rho''(1)$ we proceed in a similar way,
$$
\rho''(y) = \frac{\partial^2\Psi}{\partial x\partial
y}(\tau(y),y)\tau'(y)+ \frac{\partial^2\Psi}{\partial
y^2}(\tau(y),y).
$$
Computing the partial derivatives of $\Psi$ poses no problem; to
obtain $\tau'(y)$ we differentiate with respect to $y$ the
equation
$$
 T(\tau(y), D(\tau(y), y))=0,
$$
where $T(u,D)$ is the polynomial of Lemma~\ref{eq-def-tau}. This
gives a linear equation in $\tau'(y)$, from where it follows that
$\tau'(1)\approx -0.21992$ and then  $\rho''(1)\approx 0.44298$.
Finally, the constants $\kappa$ and $\lambda$ are computed using
$$
  \kappa = -{\rho'(1) \over \rho(1)}, \qquad
  \lambda = -{\rho''(1) \over \rho(1)}  -{\rho'(1) \over \rho(1)}
    + \left( {\rho'(1) \over \rho(1)}\right)^2.
   $$
\end{proof}

\section{Outerplanar graphs}\label{sec:outer}

We keep the notations of previous sections but applied to
outerplanar graphs instead of series-parallel graphs. Thus $g_n$
is the number of (labelled) outerplanar graphs on $n$ vertices;
similarly for $c_n$ and $b_n$, and for the corresponding
generating functions. Our exposition will be brief since the
necessary machinery has been introduced in the previous section
and the generating functions in this case are much simpler.
Moreover, we restrict our computations to the most relevant
issues, namely, the asymptotic expressions for the number of
outerplanar graphs, and the mean and variance of the expected
number of edges in a random outerplanar graphs.

As mentioned in the introduction, a 2-connected outerplanar can be
seen as a dissection of a convex polygon. The ordinary GF $A(x,y)$
for polygon dissections, where $x$ marks vertices and $y$ edges
can be easily obtained with the method introduced in \cite{FN} for
counting polygon dissections with respect to the number of faces.
Indeed, let $K$ be a convex polygon with $n$ vertices and fix an
edge $e$ of $K$. A~dissection of $K$ is either a single edge or an
ordered sequence of $k \ge2$ dissections along a face containing
$e$, where $k-1$ pairs of vertices are identified.
 Thus we have
 $$
 A(x,y) = yx^2 + y \sum_{k\ge2} {A(x,y)^k \over x^{k-1}} =
 yx^2 + {y A^2 \over x-A}.
 $$
The solution to the previous equation with non-negative terms is
 $$
    A(x,y) = \sum a_{n,k} y^kx^n = {x(1+yx - \sqrt{1-2yx - 4y^2x + y^2x^2})
    \over 2+2y}.
 $$
Returning to outerplanar graphs, each dissection of $K$ gives
raise to $(n-1)!/2$ (the number of labellings of a non-oriented
cycle) 2-connected outerplanar graphs, except for the special case
$n=2$. Hence
 $$
 b_{n,q} = a_{n,q} (n-1)! /2, \quad n\ge3, \qquad
 \hbox{and} \qquad b_{2,1} = 1.
 $$
In terms of the corresponding generating functions (recall that
$B(x,y)$ is an exponential GF), the former relations translate
into
 $$
    B'(x,y) = {{1 \over 2}A(x,y) + yx \over 2}
     = {1 + xy(3+2y) - \sqrt{1-2xy-4y^2x+x^2y^2} \over 4(1+y).}
    $$
For $y=1$, the smallest positive root of the radicand $1-6x+x^2$
is
$$R=3-2\sqrt 2 \approx 0.17157,$$
 which is then the
radius of convergence of $B(x) = B(x,1)$.

Since a graph is outerplanar if and only if its connected
components are outerplanar, and the blocks in the components are
also outerplanar, the relations we had in the previous section
between $B,C$ and $G$ also hold, that is
 $$
  G(x) = \exp(C(x)), \qquad xC'(x) =
    x\exp\left(B'(xC'(x))\right).
$$
 It follows that, for fixed $y$, the functional inverse of $F(x,y)=x C'(x,y)$ is
$$
  \Psi_y(u) = u e^{-B'(u,y)}.
$$
Given the explicit form we have for $B'$, it is easy to check that
$\Psi_1'(u)$ has a root $\tau = 0.17076$ in $[0,R]$; in fact,
$\tau$ is the smallest positive root of the equation
 $$
 3u^4 - 28u^3 + 70u^2 - 58u + 8 = 0.
 $$
Consequently, the radius of convergence of $F(x,1)$ and $C(x)$ is
equal to
 $$
    \rho = \Psi(\tau) = 0.13659.
 $$

\begin{theorem}
The number $h_n$ of outerplanar graphs is asymptotically
$$
 h_n \sim h\cdot n^{-5/2}\cdot \rho^{-n} n!
$$
where $\rho\approx 0.13659$ and $h \approx 0.017657$.
\end{theorem}

\begin{proof}
The value of the dominant singularity $\rho$ has been determined
previously. Since $\Psi'$ has a root in its domain of definition,
the inverse function $F(x)=xC'(x)$ has a singular expansion of
square-root type in $X=\sqrt{1-x\rho}$, hence $C(x)$ has  a
singular expansion whose dominant term is $X^{3/2}$ and from this
follows the subexponential term $n^{-5/2}$. Finally, the constant
$h$ is computed as in the previous section using the evaluation of
$\Psi''(\tau)$.
\end{proof}

The proof of the next theorem is omitted, since it follows exactly
the same pattern as the proof of Theorem \ref{th:edges}.
Computations in this case are simpler due to the explicit
expression for $\Psi$.

\begin{theorem}\label{th:edges-out}
Let $X_n$ denote the number of edges in a random outerplanar graph
with $n$ vertices. Then $X_n$ is asymptotically normal and the
mean $\mu_n$ and variance $\sigma_n^2$ satisfy
\begin{equation}\label{mean-var-out}
\mu_n \sim \zeta n, \qquad \sigma_n^2 \sim \eta n,
\end{equation}
where $\zeta \approx 1.56251$ and $\eta \approx 0.22399$. The same
is true, with the same constants, for \emph{connected} random
outerplanar graphs.
\end{theorem}

\section{The number of connected components}\label{sec:comp}

In this section we determine limit laws for the number of
connected components. A sequence $X_n$ of discrete random
variables converges to a discrete random variable $X$ if, for
every integer $k$,
 $$
 {\rm Prob}\{ X_n = k \} \to {\rm Prob}\{ X = k\}, \quad \hbox{ as
 $n \to \infty$}.
 $$
In the next statement, convergence is to a \emph{shifted} Poisson
law because the number of components is always strictly positive.

\begin{theorem}
The distribution of the number of connected components in random
series-parallel graphs is asymptotically  a shifted Poisson law
$1+P(\nu)$ with parameter equal to $\nu \approx 0.11761$. The same
result holds for outerplanar graphs, in this case the parameter of
the Poisson law being equal to $\xi \approx 0.14840$. As a
consequence the probability that a random series-parallel graph is
connected tends to $e^{-\nu} \approx 0.88904$, and to $e^{-\xi}
\approx 0.86208$ for outerplanar graphs.
\end{theorem}

\begin{proof}
We follow the same approach as the proof of Theorem~6 in
\cite{GN}. We present the proof for SP graphs, the case of
outerplanar graphs is analogous.

For fixed $k$, the generating function of SP graphs with exactly
$k$ connected components is
 $$
    {1 \over k!} C(x)^k.
 $$
For fixed $k$ we have
 $$
    [x^k] C(x)^k  \sim kC_0^{k-1} [x^n] C(x).
 $$
Hence the probability that a random planar SP has exactly $k$
components is asymptotically
$$
{ [x^n] C(x) / k! \over [x^n]G(x)} \sim {k C_0^{k-1} \over k!} \,
e^{-C_0} = {\nu^{k-1} \over (k-1)!} \, e^{-\nu}.
$$
If we let $\nu = C(\rho) = C_0$, the evaluation of $C(x)$ at its
dominant singularity, then the previous expression implies
convergence to a shifted Poisson law of parameter $\nu$. Since we
know the local developments around the dominant singularity, we
can compute $C_0$ exactly.
\end{proof}

\section{Graphs without a $K_{2,3}$ minor}\label{sec:k23}

In this section we analyze briefly the class of graphs that do not
contain $K_{2,3}$ as a minor; equivalently, since $K_{2,3}$ has
maximum degree three, graphs that do not contain $K_{2,3}$ as a
subdivision. They form a class strictly larger than the class of
outerplanar graphs; as we are going to see, they are not far from
this class.

Let $G$ be a 2-connected graph not containing $K_{2,3}$ as a
minor. Then either $G$ is outerplanar (no $K_4$ minor) or else $G$
contains $K_4$ as a minor, hence also as a subdivision. But if we
subdivide just one edge of $K_4$, a $K_{2,3}$ minor shows up.
Hence this subdivision can be only $K_4$. If $G$ contains an
additional vertex $x$ then, by 2-connectivity there exist an edge
$yz$ of the given $K_4$ and two internally disjoint paths from $x$
to $y$ and $z$, hence again we have a $K_{2,3}$ minor.

In conclusion if $G$ is 2-connected and does not contain $K_{2,3}$
as a minor, then either $G$ is outerplanar or $G=K_4$. If we apply
the notation of the previous section to the new class, then the
generating function $B(x,y)$ is the same as for the class of
outerplanar graphs, plus the addition of a single monomial
$y^6x^4/4!$ corresponding to the exceptional graph $K_4$. Hence
 $$
    B'(x,y) = {y^6 x^3 \over 6} +  {1 + xy(3+2y) - \sqrt{1-2xy-4y^2x+x^2y^2} \over 4(1+y).}
    $$
From this, following the same steps as in the previous section, we
obtain the following. Details are omitted in order to avoid
repetition.

\begin{theorem}
The number $s_n$ of  graphs not containing $K_{2,3}$ as a minor is
asymptotically
$$
 s_n \sim s\cdot n^{-5/2}\cdot \rho^{-n} n!
$$
where $\rho\approx 0.13648$ and $s \approx 0.013961$.
\end{theorem}

\begin{theorem}
Let $X_n$ denote the number of edges in a random graph not
containing $K_{2,3}$ as a minor with $n$ vertices. Then $X_n$ is
asymptotically normal and the mean $\mu_n$ and variance
$\sigma_n^2$ satisfy
\begin{equation}
\mu_n \sim \zeta n, \qquad \sigma_n^2 \sim \eta n,
\end{equation}
where $\zeta \approx 1.56325$ and $\eta \approx 0.224206$. The
same is true, with the same constants, for \emph{connected} random
graphs not containing $K_{2,3}$ as a minor.
\end{theorem}


\section{Concluding remarks}

We conclude with a table showing the values of the main parameters
for the classes studied in this paper, together with the class of
planar graphs studied in~\cite{GN}. For a given class of graphs
$\mathcal{G}$, the growth constant $\gamma$ is $\rho^{-1}$, where
$\rho$ is the dominant singularity of the associated generating
function $G(x)$. In other words, $\gamma=\lim_{n\to\infty}
(g_n/n!)^{1/n}$. In the last column we display the constant
$\kappa$ such that the expected number of edges is $\sim \kappa
n$.

$$
\begin{tabular}{|l|c|c|}
\hline
Class of graphs & Growth constant & Expected number of edges \\
\hline
Planar &  \hglue -5.0pt 27.2268  &  2.2132 \\
Series prallel  &  9.0733 &  1.6167\\
Outerplanar & 7.3209 & 1.5625\\
No $K_{2,3}$ minor & 7.3270& 1.5632 \\
\hline
\end{tabular}
$$

\bigskip

\bibliographystyle{amsplain}

\begin{thebibliography}{10}

\bibitem{bender} E. A. Bender, Z. Gao, N. C. Wormald, \textit{The number of
2-connected labelled planar graphs}, Elec. J. Combinatorics 9
(2002), \#43.

\bibitem{gavoille}
N. Bonichon, C. Gavoille, N. Hanusse, \textit{Canonical
Decomposition of Outerplanar Maps and Application to Enumeration,
Coding and Generation}, Springer Lecture Notes in Computer
Science, vol.. 2280, pages 81--92, 2003.

\bibitem{FN} P. Flajolet, M. Noy, \textit{Analytic Combinatorics of Non-crossing Configurations},
Discrete Math.  (1999).

\bibitem{FO} P. Flajolet, A. Odlyzko,
\textit{Singularity analysis of generating functions}, SIAM J.
Discrete Math.  3  (1990),  216--240.

\bibitem{FS} P. Flajolet, R. Sedgewick, \textit{Analytic Combinatorics} (book in
preparation), preliminary version available at {\tt
http://algo.inria.fr/flajolet/Publications  }

\bibitem{gm} S. Gerke, C. McDiarmid, \textit{On the Number of Edges in
Random Planar Graphs}, Comb. Prob. and Computing 13 (2004),
165--183.

\bibitem{GN} O. Gim\'{e}nez, M. Noy, \textit{Asymptotic enumeration
and limit laws of planar graphs}, arXiv math.CO/0501269, 14 pages.



\bibitem{moon} J.W. Moon, \emph{Counting labelled trees}, Canadian mathematical monographs, Canadian Mathematical
Congress, Montreal, Que. 1, 1970.

\bibitem{walsh} T. R. S. Walsh,
\textit{Counting labelled three-connected and homeomorphically
irreducible two-connected graphs}, J. Combin. Theory Ser. B 32
(1982), 1--11.

\end{thebibliography}

\end{document}